\documentclass[10pt,a4paper]{article} 
\usepackage{graphicx,amsmath,bm}
\usepackage{amssymb,amscd}
\usepackage{caption2,color}

\usepackage[dvipdfm,
           pdfstartview=FitH,
           CJKbookmarks=true,
           bookmarksnumbered=true,
           bookmarksopen=true,
          colorlinks=true, 
          colorlinks=cyan,
           pdfborder=001,   
           citecolor=blue%
           ]{hyperref}

\makeatletter

\@addtoreset{equation}{section} \makeatother

\newtheorem{lemma}{Lemma}[section]

\newtheorem{theorem}{Theorem}[section]
\newcommand{\dx}{\,\mathrm{d}x}
\newcommand{\ds}{\,\mathrm{d}s}
\newcommand{\n}{\nabla}
\newcommand{\p}{\partial}
\newcommand{\norm}[1]{\lVert#1\rVert}
\newcommand{\seminorm}[1]{\lvert#1\rvert}
\newcommand{\divt}{\mathrm{div}_\Gamma}      
\newcommand{\rn}{\ensuremath{\mathbb{R}^N}}
\DeclareMathAlphabet{\mathsfsl}{OT1}{cmss}{m}{sl}

\renewcommand{\vec}[1]{\mbox{\boldmath$#1$}}
\newcommand{\oo}{\ensuremath{\Omega}}
\newcommand{\G}{\Gamma}
\newcommand{\defmath}{{\,\stackrel{\mbox{\rm\tiny def}}{=}\,}}
\newcommand{\diff}{\,\mathrm{d}}

\newcommand{\mdiv}{\,\mathrm{div}}

\newcommand{\md}{\mathrm{D}}
\newcommand{\mg}{\mathcal{G}}
\newcommand{\mt}{\mathcal{T}}
\newcommand{\mq}{\mathcal{Q}}
\newcommand{\mpp}{\mathcal{P}}
\newcommand{\my}{\mathcal{Y}}
\newcommand{\mv}{\mathcal{V}}

\newcommand{\vyt}{\vec{y}_t}
\newcommand{\vvt}{\vec{V}_t}
\newcommand{\vv}{\vec{V}}
\newcommand{\vvv}{\vec{v}}

\newcommand{\vy}{\vec{y}}
\newcommand{\vyd}{\vec{y}_d}
\newcommand{\vp}{\vec{p}}

\newcommand{\vf}{\vec{f}}
\newcommand{\vg}{\vec{g}}
\newcommand{\vphi}{\vec{\varphi}}

\newcommand{\vpsi}{\vec{\psi}}

\newcommand{\vn}{\vec{n}}
\begin{document}

\title{Optimal Shape Design for Stokes Flow \\Via Minimax Differentiability\footnote{This work was
supported by the National Natural Science Fund of China under grant
number 10371096 for ZM Gao and YC Ma.}}

\author{Zhiming Gao\thanks{School of Science, Xi'an Jiaotong University, Shaanxi, P.R.China, 710049.
E--mail\,:\, dtgaozm@gmail.com.}\qquad
 Yichen Ma\footnote{Corresponding author. School of Science, Xi'an Jiaotong University, Shaanxi, P.R.China, 710049. E-mail:\,ycma@mail.xjtu.edu.cn.}
 \qquad Hongwei Zhuang\thanks{Engineering College of Armed Police Force, Shaanxi,\,P.R.China, 710086.}}
\date{}
 \maketitle
\noindent{{\textbf{Abstract.\;}} This paper is concerned with a
shape sensitivity analysis of a viscous incompressible fluid driven
by Stokes equations with nonhomogeneous boundary condition. The
structure of shape gradient with respect to the shape of the
variable domain for a given cost function is established by using
the differentiability of a minimax formulation involving a
Lagrangian functional combining with function space parametrization
technique or function space embedding technique. We apply an
gradient type algorithm to our problem. Numerical examples show that
our theory is useful for
practical purpose and the proposed algorithm is feasible. \\[8pt]
{{\textbf{Keywords.\;}}
shape optimization; minimax formulation; gradient algorithm; Stokes equations.\\[8pt]
{{\textbf{AMS(2000) subject classifications.\;}}49J35, 49K35, 49K40,
35B37.
\section{Introduction}
This paper deals with the optimal shape design for Stokes flow
inside a moving domain. This problem is a basic tool in the design
and control of many industrial devices such as aircraft wings,
automobile shapes, boats, and so on. The control variable is the
shape of the domain, the object is to minimize a cost function that
may be given by the designer, and finally we can obtain the optimal
shapes.

The efficient computation of optimal shapes requires a {shape
calculus} (see \cite{delfour}) which differs from its analog in
vector spaces. It is necessary to make sense of {shape gradient}
which is a basic tool to obtain necessary conditions and to provide
us with gradient information required by the gradient type
optimization methods. The {velocity method} (see J.Cea\cite{ce81}
and J.-P.Zolesio\cite{delfour,zo79}) gave a precise mathematical
meaning to this notion.

Many shape optimization problems can be expressed as a minimax of
some suitable Lagrangian functional. The characterization of the
change in geometric domain is obtained by velocity method. Finally
the use of theorems on the differentiability of a saddle point
(i.e., a minimax) of such lagrangian functional with respect to a
parameter provides very powerful tools to obtain {shape gradient} by
{function space parametrization} or {function space embedding}
(see\cite{delfour88}) without the usual study of the derivative of
the state.

The function space parametrization technique and function space
embedding technique are advocated by M.C.Delfour and
J.-P.Zol\'{e}sio to solving poisson equation with Dirichlet and
Nuemann condition (see\cite{delfour}). In our paper \cite{gao06}, we
apply them to a Robin problem and give its numerical implementation.
The purpose of this paper is to use lagrangian formulation and
theorem on the differentiability of a minimax to study the shape
sensitivity analysis for Stokes flow, and then give a gradient type
algorithm with some numerical examples to prove that our theory
could be very useful for the practical purpose.

This paper is organized as follows. Section 2 is devoted to the
statement of a shape optimization problem for Stokes flow. In
section 3, we briefly recall the velocity method which is used for
the characterization of the deformation of the shape of the domain,
and we also give the definitions of {Eulerian derivative} and {shape
gradient}. Then we include the divergence free condition directly
into the Lagrange functional thanks to a multiplier which plays the
role of the adjoint state associated with the primal pressure. This
leads to a saddle point formulation of the shape optimization
problem for Stokes equations with nonhomogeneous boundary condition.

 Section 4 is devoted to the computation of the
shape gradient of the Lagrangian functional due to a minimax
principle concerning the differentiability of the minimax
formulation(see\cite{correa,delfour88}) by {Function Space
Parametrization} technique.

In section 5, we compute the shape gradient by using such minimax
principle coupling with {Function Space Embedding} technique and get
the same expression obtained in section 4.

Finally, in the last section, with the shape gradient information,
we can establish a gradient type algorithm to solve our problem, and
numerical examples show the feasibility of our approach for
different viscosity coefficients.

Before closing this section, we introduce some notations that will
be used throughout the paper.

$H^m(D), m\in\mathbb{R},$ denotes the standard Sobolev space of
order $m$ with respect to the set $D$, where $D$ is either the fluid
domain $\oo$ or its boundary $\G$. Note that $H^0(D)=L^2(D)$.
Corresponding Sobolev spaces of vector-valued functions will be
denoted by $H^m(D)^N$.

Let $\vec u=(u_1, u_2,\cdots,u_d)$ and $\vec v=(v_1,
v_2,\cdots,v_d)$ be two vector functions of dimension $d$, and $w$
be a scalar function. $\md\vec u$ denotes the Jacobian matrix of
$\vec u$, i.e., $\md\vec u\defmath (\p_j u_i)_{i,j=1}^d$, and its
transpose matrix is denoted by ${}^*\md\vec u$. We also have the
following linear forms:
\begin{eqnarray*}
  &(\vec u,\vec v)_\oo\defmath\int_\oo\vec u\cdot\vec v\dx=
  \int_\oo\langle\vec u,\vec v\rangle\dx=\int_\oo\sum
  \limits^d_{i=1}u_i\, v_i\dx,\quad\forall\vec u,\vec v\in L^2(\oo)^d;\\
   &a(\oo;\vec u,\vec v)\defmath\int_{\oo}\alpha\,\md\vec u:\md\vec
  v\dx=\int_{\oo}\alpha\sum\limits^d_{i,j=1}\p_j u_i\,\p_j
  v_i\dx,\quad\forall \vec u,\vec v\in H^1(\oo)^d;\\
  &b(\oo;\vec u,w)\defmath-\int_{\oo}\mdiv\vec u\,w\dx=-\int_\oo\sum\limits_{i=1}^d\p_i u_i\,w\dx,\quad\forall\vec u\in H^1(\oo)^d,\;\forall w\in L^2(\oo).
\end{eqnarray*}
Note that the inner products in $L^2(\oo)^d$ is denoted by
$(\cdot,\cdot)_\oo$, and the angle product
$\langle\cdot,\cdot\rangle$ denotes the usual dot product of two
vectors in this paper.

\section{Formulation of the problem}
Let $\oo$ be the fluid domain in \rn ($N=2\;\mbox{or}\;3$), and the
boundary $\G\defmath\p\oo$ be smooth. The fluid is described by its
velocity $\vy$ and pressure $p$ satisfying the Stokes equations:
\begin{equation}
  \left\{%
  \begin{array}{ll}
-\alpha\Delta\vy+\n p=\vf&\quad\mbox{in}\;\oo\\
\mdiv \vy=0&\quad\mbox{in}\;\oo\\
\vy=\vg&\quad\mbox{on}\;\G
  \end{array}%
  \right.
\end{equation}
where $\alpha$ stands for the kinematic viscosity coefficient. Let
$\vf\in {H^m(\rn)}^N,$ and $\vg\in {H^{m+3/2}(\rn)}^N$ ($m\geq 0$ to
be specified) be given satisfying the compatibility condition
\begin{equation}\label{comp}
  \int_\G\vg\cdot\vn\ds=0,
\end{equation}
then we know that the solution $(\vy,p)$
belongs to $H^1(\oo)^N\times L^2(\oo)$ and even to
$H^{m+2}(\oo)^N\times H^{m+1}(\oo)$ when $\G$ is of class $C^{m+2}$
by the regularity theorem (see\cite{gilbarg,temam01}).

Our objective is to compute the first order "derivative" of the cost
function
\begin{equation}\label{prob:cost}
  J(\oo)=\frac{1}{2}\int_{\oo}\seminorm{\vy(\oo)-\vyd}^2\dx
\end{equation}
with respect to the variational domain $\oo$. The target velocity
$\vyd$ is fixed in ${H^1(\rn)}^N$ and given by the designer for some
purposes.

\section{The velocity method and a saddle point formulation}
Domains $\oo$ don't belong to a vector space and this requires the
development of {shape calculus} to make sense of a ``derivative" or
a ``gradient". To realize it, there are about three types of
techniques: J.Hadamard\cite{ha07}'s normal variation method, the
{perturbation of the identity} method by J.Simon\cite{si80} and the
{velocity method}(see J.Cea\cite{ce81} and
J.-P.Zolesio\cite{delfour,zo79}). We will use the velocity method
which contains the others.

Let $\vec V\in \mathrm{E}^k :=C([0,\tau);\mathcal{D}^k(\rn,\rn))$,
where $\mathcal{D}^k(\rn,\rn)$ denotes the space of all $k-$times
continuous differentiable functions with compact support contained
in $\rn$ and $\tau$ is a small positive real number. The velocity
field
$$\vec V(t)(x)=\vec V(t,x), \qquad x\in \rn,\quad t\geq 0$$
belongs to $\mathcal{D}^k(\rn,\rn)$ for each $t$. It can generate
transformations
 $$T_t(\vec V)X=x(t,X),\quad t\geq 0,\quad X\in\rn$$
through the following dynamical system
\begin{equation}\label{dynamical}
  \left\{%
  \begin{array}{ll}
  \frac{\diff x}{\diff t}(t,X)=\vec V(t,x(t))\\[3pt]
  x(0,X)=X
  \end{array}%
  \right.
\end{equation}
with the initial value $X$ given. We denote the "transformed domain"
$T_t(\vec V)(\oo)$ by $\oo_t(\vec V)$ at $t\geq 0$.

 Furthermore, for sufficiently small $t>0,$ the Jacobian $J_t$ is
 strictly positive:
 \begin{equation}\label{jacobian}
   J_t(x):=\det\seminorm{\md T_t(x)}=\det\md T_t(x)>0,
 \end{equation}
where $\md T_t(x)$ denotes the Jacobian matrix of the transformation
$T_t$ evaluated at a point $x\in\rn$ associated with the velocity
field $\vec V$. We will also use the following notation: $\md
T_t^{-1}(x)$ is the inverse of the matrix $\md T_t(x)$ , ${}^*\md
T_t^{-1}(x)$ is the transpose of matrix $\md T_t^{-1}(x)$, and the
Jacobian matrix of $T_t$ with respect to the boundary $\G$ is
denoted by $w_t=J_t\seminorm{{}^*\md T_t^{-1}\,\vn}$.

We now consider the solution $(\vyt,p_t)$ on $\oo_t$ of the problem
\begin{equation}\label{vm:bvpt}
\left\{%
  \begin{array}{ll}
    -\alpha\Delta \vyt+\n p_t=\vf\qquad&\mbox{in}\;\oo_t\\
    \mdiv \vyt=0,\qquad&\mbox{in}\;\oo_t\\
    \vyt=\vg\qquad&\mbox{on}\;\G_t\;(\mbox{the boundary of}\;\oo_t)
  \end{array}%
  \right.
\end{equation}
and the associated cost function
\begin{equation}\label{vm:cost}
  J(\oo_t)=\frac{1}{2}\int_{\oo_t}\seminorm{\vyt-\vyd}^2\dx
\end{equation}
We say that this functional has a \textit{Eulerian derivative} at
$\oo$ in the direction $\vec V$ if the limit
\begin{equation*}
\lim_{t\searrow 0}\frac{J(\oo_t)-J(\oo)}{t}\defmath\diff J(\oo;\vec
V)
\end{equation*}
exists.

 Furthermore, if the map
 $$\vec V\mapsto\diff
J(\oo;\vec V):\;\mathrm{E}^k\rightarrow\mathbb{R}$$ is linear and
continuous, we say that $J$ is \textit{shape differentiable} at
$\oo$. In the distributional sense we have
\begin{equation}
  \diff J(\oo;\vec V)=\langle \Im,\vec V\rangle_{\mathcal{D}^k(\rn,\rn)'\times \mathcal{D}^k(\rn,\rn)}.
\end{equation}
 When $J$ has a Eulerian derivative, we say that $\Im$ is the \textit{shape gradient} of $J$
at $\oo$.

Now we shall describe how to build an appropriate Lagrangian
functional that takes into account the divergence condition and the
nonhomogeneous Dirichlet boundary condition.

Given $\vf\in {H^1(\rn)}^N$ and $\vg\in {H^{5/2}(\rn)}^N$,  we
introduce a Lagrange multiplier $\vec\mu$ and a functional
\begin{equation}
  L(\oo,\vy,p,\vvv,q,\vec\mu)=\int_\oo \langle\alpha\Delta\vy-\n p+\vf,\vvv
  \rangle\dx+\int_\oo\mdiv\vy q\dx+\int_\G\langle\vy-\vg,\vec\mu\rangle\ds
\end{equation}for
$(\vy,p)\in Y(\oo)\times Q(\oo)$, $(\vvv,q)\in P(\oo)\times Q(\oo)$,
 and $\vec\mu\in
H^{-1/2}(\G)^N$ with
\begin{equation*}
Y(\oo_t)\defmath H^2(\oo_t)^N; \qquad P(\oo_t)\defmath
H^2(\oo_t)^N\cap H^1_0(\oo_t)^N;\quad Q(\oo_t)\defmath H^1(\oo_t),
\end{equation*}
and $\oo_0=\oo$ as $t=0$.

Now we're interested in the following saddle point problem
\begin{equation*}
  \inf_{(\vy,p)\in Y(\oo)\times Q(\oo)}\quad\sup_{(\vvv,q,\vec\mu)\in P(\oo)\times Q(\oo)\times H^{-1/2}(\G)^N}\;L(\oo,\vy,p,\vvv,q,\vec\mu)
\end{equation*}
The solution is characterized by the following:

$\bullet$\; The state $(\vy,p)$ is the solution of problem
\begin{equation}
  \left\{%
  \begin{array}{ll}
-\alpha\Delta\vy+\n p=\vf&\quad\mbox{in}\;\oo\\
\mdiv \vy=0&\quad\mbox{in}\;\oo\\
\vy=\vg&\quad\mbox{on}\;\G
  \end{array}%
  \right.
\end{equation}

$\bullet$\; The adjoint state $(\vvv,q)$ is the solution of problem
\begin{equation}\label{prob:bvp1}
\left\{%
  \begin{array}{ll}
    -\alpha\Delta \vvv+\n q=0\qquad&\mbox{in}\;\oo\\
\mdiv \vvv=0&\mbox{in}\;\oo\\
    \vvv=0\qquad&\mbox{on}\;\G;
  \end{array}%
  \right.
\end{equation}

 $\bullet$\; The multiplier satisfies: $\vec\mu=\alpha\md\vvv\,\vec{n}-q\,\vn,\;\mbox{on}\;\G$.

Hence we obtain the following new functional,
\begin{equation*}
  L(\oo,\vy,p,\vvv,q)=\int_\oo \langle\alpha\Delta\vy-\n p+\vf,\vvv
  \rangle\dx+\int_\oo\mdiv\vy q\dx+\int_\G\langle\vy-\vg,
  \alpha\md\vvv\,\vec{n}-q\,\vn\rangle\ds.
\end{equation*}
To get rid of the boundary integral, the following identities are
derived by Green formula,
\begin{eqnarray*}
\int_\G\langle\vy-\vg,\md\vvv\vn\rangle\ds&=&\int_\oo[\langle\vy-\vg,\Delta\vvv\rangle+\md(\vy-\vg):\md\vvv]\dx;\\
  \int_\G\langle\vy-\vg, q\,\vec{n}\rangle\ds&=&\int_\oo
  [\mdiv(\vy-\vg)q+\langle\vy-\vg,\n q\rangle]\dx.
\end{eqnarray*}
Thus we obtain the new Lagrangian:
\begin{multline*}
  L(\oo,\vy,p,\vvv,q)=\int_\oo \langle\alpha\Delta\vy-\n p+\vf,\vvv
  \rangle\dx+
  \alpha\int_\oo[\langle\vy-\vg,\Delta\vvv\rangle+\md(\vy-\vg):\md\vvv]\dx\\
  +\int_\oo\mdiv\vy q\dx-\int_\oo
  [\mdiv(\vy-\vg)q+\langle\vy-\vg,\n q\rangle]\dx.
\end{multline*}
 This domain integral is advantageous for the
computation of shape gradient.

Given a velocity field $\vec V\in \mathrm{E}^1$ and transformed
domain $\oo_t$, we can easily verify
\begin{equation}\label{lf:jt}
  J(\oo_t)= \inf_{(\vy_t,p_t)\in Y(\oo_t)\times Q(\oo_t)}\quad\sup_{(\vvv_t,q_t)\in P(\oo_t)\times Q(\oo_t)}G(\oo_t,\vy_t,p_t,\vvv_t,q_t)
\end{equation}
where the Lagrangian is given by
\begin{eqnarray*}
G(\oo_t,\vy_t,p_t,\vvv_t,q_t)&=&F(\oo_t,\vy_t)+L(\oo_t,\vy_t,p_t,\vvv_t,q_t)\\
                    &=&\frac{1}{2}\int_{\oo_t}\seminorm{\vy_t-\vyd}^2\dx
                    +\int_{\oo_t}\langle\alpha\Delta\vy_t-\n p_t+\vf,\vvv_t
  \rangle\dx+\int_{\oo_t}\mdiv\vy_t\, q_t\dx\\
  &&+
  \alpha\int_{\oo_t}[\langle\vy_t-\vg,\Delta\vvv_t\rangle+\md(\vy_t-\vg):\md\vvv_t]\dx\\
  &&-\int_{\oo_t}
  [\mdiv(\vy_t-\vg)q_t+\langle\vy_t-\vg,\n q_t\rangle]\dx.
\end{eqnarray*}
and $J(\oo_t)$ was characterized by (\ref{vm:cost}).

The Lagrangian $G(\oo_t,\cdot,\cdot,\cdot,\cdot)$ has a unique
saddle point $(\vyt,p_t,\vvv_t,q_t)\in Y(\oo_t)\times Q(\oo_t)\times
P(\oo_t)\times Q(\oo_t)$ which is given by the following systems:
\begin{description}
    \item[State equations]
    \begin{subequations}\label{vm1}
\begin{equation}\label{vm:bvpu}
  \left\{%
  \begin{array}{ll}
-\alpha\Delta\vy_t+\n p_t=\vf&\quad\mbox{in}\;\oo_t\\
\mdiv \vy_t=0&\quad\mbox{in}\;\oo_t\\
\vyt=\vg&\quad\mbox{on}\;\G_t
  \end{array}%
  \right.
\end{equation}
    \item[Adjoint state equations]
\begin{equation}\label{vm:bvpp}
\left\{%
  \begin{array}{ll}
    -\alpha\Delta \vvv_t+\n q_t=\vy_t-\vyd\qquad&\mbox{in}\;\oo_t\\
\mdiv \vvv_t=0&\mbox{in}\;\oo_t\\
    \vvv_t=0\qquad&\mbox{on}\;\G_t;
  \end{array}%
  \right.
\end{equation}
\end{subequations}
\end{description}
Our objective is to get the limit
\begin{equation}
  \diff j(0)=\lim_{t\searrow 0}\frac{j(t)-j(0)}{t}
\end{equation}
where $j(t)=J(\oo_t)= \inf_{(\vy_t,p_t)\in Y(\oo_t)\times
Q(\oo_t)}\quad\sup_{(\vvv_t,q_t)\in P(\oo_t)\times
Q(\oo_t)}G(\oo_t,\vy_t,p_t,\vvv_t,q_t)$.

Unfortunately, the Sobolev space $Y(\oo_t)$, $Q(\oo_t)$, and
$P(\oo_t)$ depend on the parameter $t$, so we need a theorem to
differentiate a saddle point with respect to the parameter $t$, and
there are two techniques to get rid of it:
\begin{itemize}
    \item   \textit{Function space parametrization }technique;
    \item \textit{ Function space embedding }technique.
\end{itemize}
In section \ref{fsp} we will use the first case, and section
\ref{fse} is devoted to the second case. We will find that both of
them can derive the same expression for $\diff J(\oo;\vec V)$.
\section{Function space parametrization}\label{fsp}
This section is devoted to the \textit{function space
parametrization}, which consists in transporting the different
quantities (such as, a cost function) defined on the variable domain
$\oo_t$ back into the reference domain $\oo$ which does not depend
on the perturbation parameter $t$. Thus we can use differential
calculus since the functionals involved are defined in a fixed
domain $\oo$ with respect to the parameter $t$.

We parameterize the functions in $H^m(\oo_t)^d$ by elements of
$H^m(\oo)^d$ through the transformation:
\begin{equation*}
  \vphi\mapsto \vphi\circ T_t^{-1}:\quad H^m(\oo)^d\rightarrow
  H^m(\oo_t)^d,\qquad \mbox{integer}\;m\geq 0.
\end{equation*}
where "$\circ$" denotes the composition of the two maps and $d$ is
the dimension of the function $\vphi$.

 Note that
since $T_t$ and $T_t^{-1}$ are diffeomorphisms, it transforms the
reference domain $\oo$ (respectively, the boundary $\G$) into the
new domain $\oo_t$ (respectively, the boundary $\G_t$ of $\oo_t$).
This parametrization can not change the value of the saddle point.
We can rewrite (\ref{lf:jt}) as
\begin{equation}\label{fsp:newsaddlep}
J(\oo_t)= \inf_{(\vy,p)\in Y(\oo)\times
Q(\oo)}\quad\sup_{(\vvv,q)\in P(\oo)\times Q(\oo)}G(\oo_t,\vy\circ
T_t^{-1},p\circ T_t^{-1},\vvv\circ T_t^{-1},q\circ T_t^{-1}).
\end{equation}
It amounts to introducing the new Lagrangian for $(\vy,p,\vvv,q)\in
Y(\oo)\times Q(\oo)\times P(\oo)\times Q(\oo)$:
\begin{equation*}
  \tilde G(t,\vy,p,\vvv,q)\defmath G(\oo_t,\vy\circ
T_t^{-1},p\circ T_t^{-1},\vvv\circ T_t^{-1},q\circ T_t^{-1}).
\end{equation*}
The expression for $\tilde G(t,\vy,p,\vvv,q)$ is given by
\begin{equation}\label{gtp}
  \tilde G(t,\vy,p,\vvv,q)=I_1(t)+I_2(t)+I_3(t)+I_4(t),
\end{equation} where
\begin{eqnarray*}
 I_1(t)&\defmath&\frac{1}{2}\int_{\oo_t}\seminorm{\vy\circ
 T_t^{-1}-\vyd}^2\dx;\\
 I_2(t)&\defmath&\int_{\oo_t} \langle\alpha\Delta(\vy\circ
T_t^{-1})-\n (p\circ T_t^{-1})+\vf,\vvv\circ T_t^{-1}
  \rangle\dx+\int_{\oo_t}\mdiv(\vy\circ
T_t^{-1})(q\circ T_t^{-1})\dx;\\
I_3(t)&\defmath&\alpha\int_{\oo_t}[\langle\vy\circ
T_t^{-1}-\vg,\Delta(\vvv\circ T_t^{-1})\rangle+\md(\vy\circ
T_t^{-1}-\vg):\md(\vvv\circ
T_t^{-1})]\dx;\\
I_4(t)&\defmath&-\int_{\oo_t}[\mdiv(\vy\circ T_t^{-1}-\vg)(q\circ
T_t^{-1})+\langle\vy\circ T_t^{-1}-\vg,\n(q\circ
T_t^{-1})\rangle]\dx,
\end{eqnarray*}
and its saddle point is the solution of the following variational
systems:
\begin{description}
    \item[State system] $\qquad(\vy^t,p^t)\in Y(\oo)\times Q(\oo),\qquad \forall (\vpsi,\pi)\in P(\oo)\times Q(\oo),$
    \begin{subequations}\label{ytt1}
\begin{equation}\label{ytt}
\left\{
\begin{array}{ll}
   a(\oo_t;\vy^t\circ T_t^{-1},\vpsi\circ
  T_t^{-1})+b(\oo_t;\vpsi\circ
  T_t^{-1},p^t\circ
  T_t^{-1}) =(\vf,\vpsi\circ T_t^{-1})_{\oo_t};\\
b(\oo_t;\vy^t\circ
  T_t^{-1},\pi\circ
  T_t^{-1})=0.
\end{array}
\right.
\end{equation}
    \item[Adjoint state system]$\qquad(\vvv^t,q^t)\in P(\oo)\times Q(\oo),\qquad \forall (\vphi,r)\in P(\oo)\times Q(\oo),$
    \begin{equation}\label{ptt}
    \left\{
    \begin{array}{ll}
a(\oo_t;\vvv^t\circ T_t^{-1},\vphi\circ T_t^{-1})+b(\oo_t;\vphi\circ
  T_t^{-1},q^t\circ
  T_t^{-1})=(\vy^t\circ T_t^{-1}-\vyd,\vphi\circ
T_t^{-1})_{\oo_t},\;\\
b(\oo_t;\vvv^t\circ
  T_t^{-1},r\circ
  T_t^{-1})=0.
    \end{array}
\right.
    \end{equation}
    \end{subequations}
\end{description}
By Green formula, the equivalent expression for $\tilde
G(t,\vy^t,p^t,\vvv^t,q^t)$ is obtained:
\begin{multline}\label{gtvphi}
\tilde
G(t,\vy^t,p^t,\vvv^t,q^t)=\frac{1}{2}\int_{\oo_t}\seminorm{\vy^t\circ
T_t^{-1}-\vyd}^2\dx\\
-a(\oo_t;\vy^t\circ T_t^{-1},\vvv^t\circ
T_t^{-1})-b(\oo_t;\vvv^t\circ T_t^{-1},p^t\circ
T_t^{-1})+(\vf,\vvv^t\circ T_t^{-1})_{\oo_t}\\
-b(\oo_t;\vy^t\circ T_t^{-1},q^t\circ T_t^{-1})
+\int_{\G_t}\left\langle\vy^t\circ
T_t^{-1}-\vg,\alpha\md(\vvv^t\circ T_t^{-1})\vn-(q^t\circ
T_t^{-1})\vn\right \rangle\,\ds.
\end{multline}
 By the transformation $T_t$, and the following two chain rule identities,
\begin{eqnarray*}
  \md(\vphi\circ T_t^{-1})=(\md\vphi\cdot[\md T_t]^{-1})\circ
  T_t^{-1};\\
  \mdiv (\vphi\circ T_t^{-1})=(\md\vphi:{}^*[\md T_t]^{-1})\circ
  T_t^{-1},
\end{eqnarray*}
we can rewrite it on $\oo$ as
\begin{multline}
\tilde G(t,\vy^t,p^t,\vvv^t,q^t)=\frac{1}{2}\int_\oo
\seminorm{\vy^t-\vyd\circ T_t}^2 J_t\dx\\-\int_\oo
\mathcal{A}(t)\md\vy^t:\md\vvv^t\,J_t\dx-\int_\oo \mathcal{B}(t)\md
 \vvv^t p^t\,J_t\dx+\int_\oo\langle\vf\circ
T_t,\vvv^t\rangle J_t\dx\\
-\int_\oo \mathcal{B}(t)\md
 \vy^t q^t\,J_t\dx+\int_\G \left\langle\vy^t-\vg\circ T_t,
 \mathcal{C}(t)\md \vvv^t-q^t(\vn\circ T_t)\right
\rangle\, w_t\ds
\end{multline}
where the notation
\begin{eqnarray*} \mathcal{A}(t)\vec\tau:\vec\sigma&\defmath&
\alpha[\vec\tau(\md T_t)^{-1}]:[\vec\sigma(\md T_t)^{-1}],\;\quad
  \mathcal{B}(t)\vec\tau\defmath[\vec\tau:{}^*(\md T_t)^{-1}];\\
\mathcal{C}(t)\vec\tau&\defmath&\alpha\vec\tau(\md
T_t)^{-1}(\vn\circ T_t),\qquad w_t\defmath J_t\seminorm{{}^*\md
T_t^{-1}\,\vn}.
\end{eqnarray*}
Similarly, the variational systems (\ref{ytt1}) become to
\begin{description}
    \item[State system]\qquad\qquad\quad
    $(\vy^t,p^t)\in Y(\oo)\times Q(\oo),\qquad \forall (\vpsi,\pi)\in P(\oo)\times Q(\oo),$
    \begin{subequations}\label{yt1}
    \begin{equation}\label{yt}
\left\{
\begin{array}{ll}
\int_\oo
\mathcal{A}(t)\md\vy^t:\md\vpsi\,J_t\dx+\int_\oo\mathcal{B}(t)\md\vpsi
p^t\,J_t\dx=\int_\oo\langle\vf\circ T_t,\vpsi\rangle J_t\dx.\\
\int_\oo\mathcal{B}(t)\md\vy^t \pi\,J_t\dx=0.
\end{array}
\right.
    \end{equation}
    \item[Adjoint state system]$\qquad(\vvv^t,q^t)\in P(\oo)\times Q(\oo),\qquad \forall (\vphi,r)\in P(\oo)\times Q(\oo),$
        \begin{equation}\label{pt}
        \left\{
        \begin{array}{ll}
\int_\oo
\mathcal{A}(t)\md\vvv^t:\md\vphi\,J_t\dx+\int_\oo\mathcal{B}(t)\md\vphi
q^t\,J_t\dx=\int_\oo\langle\vy^t-\vyd\circ T_t,\vphi\rangle
J_t\dx,\\
\int_\oo\mathcal{B}(t)\md\vvv^t r\,J_t\dx=0.
        \end{array}
        \right.
    \end{equation}
    \end{subequations}
\end{description}
Now we introduce the theorem concerning on the differentiability of
a saddle point (or a minimax). To begin with, some notations are
given as follows.

 Define a functional
$$\mg : [0,\tau]\times X\times Y\rightarrow\mathbb{R}$$
with $\tau>0$, and $X,Y$ are the two topological spaces.

 For any
$t\in [0,\tau]$, define
$$g(t)=\inf_{x\in X}\sup_{y\in Y}\mg(t,x,y)$$
and the sets
\begin{eqnarray*}
  &X(t)=\{x^t\in X:g(t)=\sup_{y\in Y}\mg(t,x^t,y)\}\\
  &Y(t,x)=\{y^t\in Y:\mg(t,x,y^t)=\sup_{y\in Y}\mg(t,x,y)\}
\end{eqnarray*}
Similarly, we can define dual functionals
$$h(t)=\sup_{y\in Y}\inf_{x\in X}\mg(t,x,y)$$
and the corresponding sets
\begin{eqnarray*}
 & Y(t)=\{y^t\in Y:h(t)=\inf_{x\in X}\mg(t,x,y^t)\}\\
  &X(t,y)=\{x^t\in X:\mg(t,x^t,y)=\inf_{x\in X}\mg(t,x,y)\}
\end{eqnarray*}
Furthermore, we introduce the set of saddle points
$$S(t)=\{(x,y)\in X\times Y: g(t)=\mg(t,x,y)=h(t)\}$$
Now we can introduce the following theorem (see \cite{correa} or
page 427 of \cite{delfour}):
\begin{theorem}\label{fsp:correa}
 Assume that the following hypothesis hold:
 \begin{itemize}
    \item [(H1)]$S(t)\neq\emptyset,\;t\in [0,\tau];$
    \item [(H2)]The partial derivative $\p_t\mg(t,x,y)$ exists in
    $[0,\tau]$ for all $$(x,y)\in \left[\underset{{t\in [0,\tau]}}{\bigcup}X(t)\times Y(0)\right]\bigcup\left[X(0)\times\underset{{t\in [0,\tau]}}{\bigcup}Y(t)\right];$$
    \item [(H3)]There exists a topology $\mt_X$ on $X$ such that for
    any sequence $\{t_n:t_n\in [0,\tau]\}$ with
    $\lim\limits_{n\nearrow\infty}t_n=0$, there exists $x^0\in X(0)$ and a subsequence
    $\{t_{n_k}\}$, and for each $k\geq 1,$ there exists $x_{n_k}\in
    X(t_{n_k})$ such that
    \begin{enumerate}
        \item [(i)]$\lim\limits_{n\nearrow\infty}x_{n_k}=x^0$ in the
        $\mt_X$ topology,
        \item [(ii)]$$\liminf\limits_{t\searrow 0\atop k\nearrow\infty}\p_t\mg(t,x_{n_k},y)\geq\p_t\mg(0,x^0,y),\quad \forall y\in Y(0);$$
    \end{enumerate}
    \item [(H4)]There exists a topology $\mt_Y$ on $Y$ such that for
    any sequence $\{t_n:t_n\in [0,\tau]\}$ with
    $\lim\limits_{n\nearrow\infty}t_n=0$, there exists $y^0\in Y(0)$ and a subsequence
    $\{t_{n_k}\}$, and for each $k\geq 1,$ there exists $y_{n_k}\in
    Y(t_{n_k})$ such that
    \begin{enumerate}
        \item [(i)]$\lim\limits_{n\nearrow\infty}y_{n_k}=y^0$ in the
        $\mt_Y$ topology,
        \item [(ii)]$$\limsup\limits_{t\searrow 0\atop k\nearrow\infty}
        \p_t\mg(t,x,y_{n_k})\leq\p_t\mg(0,x,y^0),\quad \forall x\in X(0).$$
    \end{enumerate}
 \end{itemize}
 Then there exists $(x^0,y^0)\in X(0)\times Y(0)$ such that
 \begin{multline}
   \diff g(0)=\lim_{t\searrow 0}\frac{g(t)-g(0)}{t}\\=\inf_{x\in
   X(0)}\sup_{y\in Y(0)}\p_t \mg(0,x,y)=\p_t\mg(0,x^0,y^0)=\sup_{y\in Y(0)}\inf_{x\in
   X(0)}\p_t \mg(0,x,y)
 \end{multline}
 This means that $(x^0,y^0)\in X(0)\times Y(0)$ is a saddle point of
 $\p_t\mg(0,x,y)$.
\end{theorem}

In order to apply Theorem \ref{fsp:correa} to our problem, we should
verify the four assumptions (H1)--(H4) below.

First of all, Let's check (H1). Assume that the velocity field $\vec
V\in \mathrm{E}^1$. Choose $\tau>0$ small enough, such that there
exists two constants $\alpha_0,\beta_0(0<\alpha_0<\beta_0)$,
$$\alpha_0\leq
  J_t(=\seminorm{J_t})\leq\beta_0,\;\forall\, t\in [0,\tau].$$
Now we can follow the standard proof of the existence and uniqueness
of solutions of Stokes equations (see \cite{temam01}) to obtain that
there exists a unique solution $(\vy^t,p^t,\vvv^t,q^t)\in
Y(\oo)\times Q(\oo)\times P(\oo)\times Q(\oo)$ ($p^t$ and $q^t$ are
unique up to a constant) and
$$\forall\,t\in [0,\tau],\quad X(t)=\{\vy^t,p^t\}\neq\emptyset,\quad Y(t)=\{\vvv^t,q^t\}\neq\emptyset.$$
Thus (H1) is satisfied.

The next step is to verify (H2). The partial derivative of $\tilde
G(\vy^t,p^t,\vvv^t,q^t)$ with respect to the parameter $t$ is
characterized by
\begin{multline}
  \p_t\tilde
G(\vy^t,p^t,\vvv^t,q^t)=\int_\oo\left
[\frac{1}{2}\seminorm{\vy^t-\vyd\circ T_t}^2 \mdiv\vec
V_t+J_t\langle \vy^t-\vyd\circ T_t,-\md \vyd\vec V_t
\rangle\right ]\dx\\
-\int_\oo [\mathcal{A}'(t)\md\vy^t:\md\vvv^t
J_t+\mathcal{A}(t)\md\vy^t:\md\vvv^t\mdiv\vec V_t]\dx\\-\int_\oo
[\mathcal{B}'(t)\md
 \vvv^t p^t\,J_t+\mathcal{B}(t)\md
 \vvv^t p^t\,\mdiv\vec V_t]\dx+\int_\oo [\langle\md\vf\vec V_t,\vvv^t\rangle J_t+\langle\vf\circ
T_t,\vvv^t\rangle\mdiv\vec
V_t]\ds\\
-\int_\oo [\mathcal{B}'(t)\md
 \vy^t q^t\,J_t+\mathcal{B}(t)\md
 \vy^t q^t\,\mdiv \vec V_t]\dx
 +\int_\G \left\langle-\md\vg\vec V_t,
 \mathcal{C}(t)\md \vvv^t-q^t(\vn\circ T_t)\right
\rangle\, w_t\dx\\
+\int_\G \left\langle\vy^t-\vg\circ T_t,
 [\mathcal{C}(t)\md \vvv^t-q^t(\vn\circ T_t)]\mdiv_\G\vec V_t+[\mathcal{C}'(t)\md \vvv^t-q^t\md\vn\vec V_t]w_t\right
\rangle\, w_t\ds.
\end{multline}
where
\begin{eqnarray*}
  \mt(t)&=&(\md T_t)^{-1},\qquad \mt'(t)=-\mt(t)\md\vec V_t\circ T_t;\\
\mathcal{B}'(t)\tau&=&\vec\tau:{}^*\mt'(t)),\qquad  \mathcal{C}'(t)\vec\tau=\alpha\vec\tau\,\mt'(t)\,(\vn\circ T_t);\\
\mathcal{A}'(t)\vec\tau:\vec\sigma&=&\alpha[(\vec\tau\mt'(t)):(\vec\sigma\mt(t))+(\vec\tau\mt(t)):(\vec\sigma\mt'(t))]\\
\divt V_t&=&\mdiv\vec V_t-\md\vec V_t\vn\cdot\vn.
\end{eqnarray*}
Since $\vec V\in \mathrm{E}^1$, $t\mapsto\vvt$ and $t\mapsto\md\vvt$
are continuous, we know that for all $(\vy^t,p^t,\vvv^t,q^t)\in
Y(\oo)\times Q(\oo)\times P(\oo)\times Q(\oo)$, $\p_t\tilde
G(\vy^t,p^t,\vvv^t,q^t)$ is well defined and exists everywhere in
$[0,\tau]$ provided that $\vf, \vyd\in H^1(\rn)^N$ and $\vg\in
H^{5/2}(\rn)^N$.

To check (H3)(i) and (H4)(i), firstly we can readily show that there
exists a positive constant $c$ such that
\begin{eqnarray*}
  \norm{\vy^t}_{H^1(\oo)^N}+\norm{p^t}_{L^2(\oo)}\leq c\norm{\vf}_{L^2(\rn)^N}\\
  \norm{\vvv^t}_{H^1(\oo)^N}+\norm{q^t}_{L^2(\oo)}\leq c\norm{\vy^t-\vyd}_{L^2(\oo)^N}
\end{eqnarray*}
Hence there exists subsequences $(\vy^{t_n},p^{t_n})$,
$(\vvv^{t_n},q^{t_n})$ and a priori $(\vec z_1,s_1)$, $(\vec
z_2,s_2)$ such that
\begin{eqnarray*}
\vy^{t_n}\rightharpoonup\vec z_1,\quad \vvv^{t_n}\rightharpoonup\vec
z_2,&&\qquad\mbox{ weakly in }H^1(\oo)^N;\\
p^{t_n}\rightharpoonup s_1,\quad q^{t_n}\rightharpoonup
s_2,&&\qquad\mbox{ weakly in }L^2(\oo).
\end{eqnarray*}
Passing to the limit, $(\vec z_1,s_1)$ is characterized by
\begin{equation*}
  \left\{
  \begin{array}{ll}
a(\oo;\vec z_1,\vpsi)+b(\oo;\vpsi,s_1)=(\vf,\vpsi)_\oo,&\quad\forall\,\vpsi\in P(\oo);\\
b(\oo;\vec z_1,\pi)=0,&\quad\forall\pi\in Q(\oo),
  \end{array}
  \right.
\end{equation*}
and $(\vec z_2,s_2)$ satisfies:
\begin{equation*}
  \left\{
  \begin{array}{ll}
a(\oo;\vec z_2,\vphi)+b(\oo;\vphi,s_2)=(\vec z_1-\vyd,\vphi)_\oo,&\quad\forall\,\vphi\in P(\oo);\\
b(\oo;\vec z_2,r)=0,&\quad\forall r\in Q(\oo),
  \end{array}
  \right.
\end{equation*}
By uniqueness, we obtain $(\vec z_1,s_1)=(\vy,p)$ and $(\vec
z_2,s_2)=(\vvv,q)$, where $(\vy,p)$ and $(\vvv,q)$ is the solution
of (\ref{ytt}) and (\ref{ptt}) at $t=0$, respectively. i.e.,
\begin{equation}\label{p}
  \left\{
  \begin{array}{ll}
a(\oo;\vec y,\vpsi)+b(\oo;\vpsi,p)=(\vf,\vpsi)_\oo,&\quad\forall\,\vpsi\in P(\oo);\\
b(\oo;\vec y,\pi)=0,&\quad\forall\pi\in Q(\oo),
  \end{array}
  \right.
\end{equation}
and
\begin{equation}\label{u}
  \left\{
  \begin{array}{ll}
a(\oo;\vec v,\vphi)+b(\oo;\vphi,q)=(\vec y-\vyd,\vphi)_\oo,&\quad\forall\,\vphi\in P(\oo);\\
b(\oo;\vec v,r)=0,&\quad\forall r\in Q(\oo).
  \end{array}
  \right.
\end{equation}
 Furthermore, we can deduce the $H^1(\oo)^N\times L^2(\oo)-$strong convergence: $
  (\vy^{t_n},p^{t_n})\rightarrow (\vy,p)$ and $(\vvv^{t_n},q^{t_n})\rightarrow
  (\vvv,q)$, Hence (H3)(i) and (H4)(i) are satisfied for the
$H^2(\oo)^N\times H^1(\oo)-$strong topology by the classical
regularity theorem(see \cite{gilbarg,temam01}). Finally, assumptions
(H3)(ii) and (H4)(ii) are readily satisfied in view of the strong
continuity of $(t,\vy,p)\mapsto \p_t\tilde G(t,\vy,p,\vvv,q)$ and
$(t,\vvv,q)\mapsto \p_t\tilde G(t,\vy,p,\vvv,q)$.

Hence all the four assumptions are satisfied, and we have the
Eulerian derivative:
\begin{multline}\label{fsp:end}
\diff J(\oo;\vec
V)=\int_\oo\left[\frac{1}{2}\seminorm{\vy-\vyd}^2\mdiv\vec
V-\langle\vy-\vyd,\md\vyd\vec V\rangle\right] \dx-\int_\oo \mdiv[(\alpha{}^*\md\vvv-q\,\mathrm{I})(\md\vg\vv)] \dx\\
-\int_\oo [\mathcal{A}'(0)\md\vy:\md\vvv+\mathcal{B}'(0)\md\vvv
p-\langle\md\vf\vec V,\vvv\rangle+\mathcal{B}'(0)\md\vy q]\dx,
\end{multline}
where $(\vy,p)$ and $(\vvv,q)$ are characterized by the variational
system(\ref{p}) and (\ref{u}), respectively, and the notation
\begin{equation*}
 \mathcal{A}'(0)\md\vy:\md\vp=-\alpha[(\md\vy\md\vec V):\md\vp+\md\vy:(\md\vp\md\vec
 V)];\quad \mathcal{B}'(0)\vec\tau=-\vec\tau:{}^*\md\vv.
\end{equation*}

Expression (\ref{fsp:end}) is a domain integral, and it is easy to
find that the map
$$\vec V\mapsto\diff J(\oo;\vec V):\quad\mathrm{E}^1\rightarrow\mathbb{R}$$
is linear and continuous, i.e., $J(\oo)$ is \textit{shape
differentiable}. Then according to Hadamard-Zol\'{e}sio structure
theorem (see \cite{delfour},Thm.3.6 and Cor.1, p.348), there exists
a scalar distribution $\mathcal{W}(\G)\in \mathcal{D}^1(\G)'$ such
that
$$\diff J(\oo;V)=\int_\G\mathcal{W}(\G)<\vec V,\vn>\ds.$$
Now we further characterize this boundary expression. Since
$(\vy,p,\vvv,q)\in H^3(\oo)^N\times H^2(\oo)\times H^3(\oo)^N\times
H^2(\oo)$ provided that $\G$ is at less $C^3$ (see \cite{temam01}),
 we can use Hadamard formula (see \cite{delfour,zolesio}):
\begin{equation}\label{hadamard}
 \frac{\diff{}}{\diff t}\int_{\oo_t}F(t,x)\dx=\int_{\oo_t}
 \frac{\p F}{\p t}(t,x)\dx+\int_{\p\oo_t} F(t,x)\langle\vec
 V,\vn_t\rangle\diff\G_t
 \end{equation}
 for a sufficiently smooth functional
$F:[0,\tau]\times\rn\rightarrow\mathbb{R}$. So we can compute the
partial derivative for $\tilde G(t,\vy,p,\vvv,q)$ with the
expression (\ref{gtp}) by using Hadamard formula,
\begin{equation*}
\frac{\p}{\p t}\left\{I_1(t)\right\}\Big{|}_{t=0}=\int_\oo
\langle\vy-\vyd,-\md\vy\vec
  V\rangle\dx+\frac{1}{2}\int_\G\seminorm{\vy-\vyd}^2<\vec
  V,\vn>\ds;
\end{equation*}
\begin{multline*}
\frac{\p}{\p t}\left\{I_2(t)\right\}\Big{|}_{t=0}=\int_{\oo}
\langle\alpha\Delta(-\md\vy\vec V)-\n(-\n p\cdot\vec
V),\vvv\rangle\dx +\int_\oo\langle\alpha\Delta\vy-\n
p+\vf,-\md\vvv\vec
V\rangle\dx\\
+\int_\oo [\mdiv(-\md\vy\vec V)q+\mdiv\vy(-\n q\cdot\vec V)]\dx\\
+\int_\G[\langle\alpha\Delta\vy-\n
p+\vf,\vvv\rangle+\mdiv\vy\,q]\langle\vec V,\vn\rangle\ds;
\end{multline*}
\begin{multline*}
\frac{\p}{\p t}\left\{I_3(t)\right\}\Big{|}_{t=0}=\alpha\int_\oo
\left\{\langle-\md\vy\vec
V,\Delta\vvv\rangle+\langle\vy-\vg,\Delta(-\md\vvv\vec
  V)\rangle\right.\\\left.{\qquad}-\md(\md\vy\vec V):\md\vvv-\md(\vy-\vg):\md(\md\vvv\vec
  V)\right\}\dx+\alpha\int_\G [\langle\vy-\vg,\Delta\vvv\rangle+\md(\vy-\vg):\md\vvv]\langle\vec
  V,\vn\rangle\ds;
\end{multline*}
\begin{multline*}
\frac{\p}{\p t}\left\{I_4(t)\right\}\Big{|}_{t=0}=-\int_\oo
[\mdiv(-\md\vy\vec V)q+\mdiv(\vy-\vg)(-\n q\cdot\vec
V)\\
{\hspace*{1cm}}-\langle\md\vy\vec V,\n q\rangle-\langle\vy-\vg,\n(\n
q\cdot\vec V)\rangle]\dx-\int_\G [\mdiv(\vy-\vg)q+\langle\vy-\vg, \n
q\rangle]\langle\vec V,\vn\rangle\ds.
\end{multline*}
Since $(\vy,p)$ and $(\vvv,q)$ are characterized by (\ref{p}) and
(\ref{u}) respectively, we obtain the boundary expression for the
shape gradient,
\begin{eqnarray}\nonumber
 \diff J(\oo;\vec V)&=&\frac{\p}{\p
t}\left\{I_1(t)+I_2(t)+I_3(t)+I_4(t)\right\}\Big{|}_{t=0}\\\label{end}
&=&\int_\G\left\{\frac{1}{2}\seminorm{\vy-\vyd}^2+\alpha\md(\vy-\vg):\md\vvv
\right\}\vec V\cdot\vn\ds.
\end{eqnarray}

\section{Function space embedding}\label{fse}
In the previous section, we have used the technique of function
space parametrization in order to get the derivative of $J(\oo_t)$,
i.e.,
\begin{equation}
J(\oo_t)=\inf_{(\vy,p)\in Y(\oo_t)\times Q(\oo_t)}\sup_{(\vvv,q)\in
  P(\oo_t)\times Q(\oo_t)}G(\oo_t,\vy,p,\vvv,q).
\end{equation}
with respect to the parameter
$t>0.$ This section is devoted to a different method based on
function space embedding technique. It means that the state and
adjoint state are defined on a large enough domain $D$ (called a
\textit{hold-all} \cite{delfour}) which contains all the
transformations $\{\oo_t: 0\leq t\leq\tau\}$ of the reference domain
$\oo$ for some small $\tau>0.$

For convenience, let $D=\rn$. Use the function space embedding,
\begin{equation}\label{fse:fun}
  J(\oo_t)=\inf_{(\vec\my,\mpp)\in Y(\rn)\times Q(\rn)}\sup_{(\vec\mv,\mq)\in
  P(\rn)\times Q(\rn)}G(\oo_t,\vec\my,\mpp,\vec\mv,\mq).
\end{equation}
where the new Lagrangian
\begin{multline}\label{fse:lag}
G(\oo_t,\vec{\mathcal Y},{\mathcal P},\vec\mv,\mathcal{Q})=F(\oo_t,\vec\my)+L(\oo_t,\vec{\mathcal Y},{\mathcal P},\vec\mv,\mathcal{Q})\\
                    =\frac{1}{2}\int_{\oo_t}\seminorm{\vec\my-\vyd}^2\dx
                    +\int_{\oo_t}\langle\alpha\Delta\vec\my-\n \mpp+\vf,\vec\mv
  \rangle\dx+\int_{\oo_t}\mdiv\vec\my\, \mq\dx\\
  +
  \alpha\int_{\oo_t}[\langle\vec\my-\vg,\Delta\vec\mv\rangle+\md(\vec\my-\vg):\md\vec\mv]\dx-\int_{\oo_t}
  [\mdiv(\vec\my-\vg)\mq+\langle\vec\my-\vg,\n \mq\rangle]\dx.
\end{multline}
Since $\vf,\vyd\in H^1(\rn)^N,$ $\vg\in H^{5/2}(\rn)^N$, and $\oo_t$
is sufficiently smooth, the unique solution $(\vy_t,p_t,\vvv_t,q_t)$
of (\ref{vm1}) belongs to $H^3(\oo_t)^N\times (H^2(\oo_t)\cap
L^2_0(\oo_t))\times (H^3(\oo_t)^N\cap H^1_0(\oo_t)^N)\times
(H^2(\oo_t)\cap L^2_0(\oo_t))$ instead of $Y(\oo_t)\times
Q(\oo_t)\times P(\oo_t)\times Q(\oo_t).$ Therefore, the sets
$$X=Y=H^3(\rn)^N\times H^2(\rn),$$
and the saddle points $S(t)=X(t)\times Y(t)$ are given by
\begin{eqnarray}
  X(t)&=&\{(\vec\my,\mpp)\in X: \,\vec\my|_{\oo_t}=\vy_t,\;\mpp|_{\oo_t}=p_t\}\\
  Y(t)&=&\{(\vec\mv,\mq)\in Y: \,\vec\mv|_{\oo_t}=\vvv_t,\;\mq|_{\oo_t}=q_t\}
\end{eqnarray}
Now we begin to verify the four assumptions of Theorem
\ref{fsp:correa} . Firstly, we can always construct a linear and
continuous extension(see \cite{adams}):
\begin{equation}
  \Pi: H^m(\oo)^d\rightarrow H^m(\rn)^d,\qquad d=1\mbox{ or }N,
\end{equation}
and
\begin{equation}
  \Pi_t: H^m(\oo_t)^d\rightarrow H^m(\rn)^d.
\end{equation}
Therefore we can define the extensions
\begin{equation}
  \vec\my_t=\Pi_t\vy_t,\;\mpp_t=\Pi_t p_t,\quad\mbox{and}\quad
  \vec\mv_t=\Pi_t\vvv_t,\;\mq_t=\Pi_t q_t,
\end{equation}
of $\vy_t$, $p_t$, $\vvv_t$ and $q_t$. So $(\vec \my_t,\mpp_t)\in
X(t)$ and $(\vec\mv_t,\mq_t)\in Y(t)$, and this shows the existence
of a saddle point, i.e.,
 $S(t)\neq\emptyset$. Then (H1) is satisfied.

To check (H2), we compute the partial derivative of the expression
(\ref{fse:lag}),
\begin{equation}\label{fse:aaa}
  \p_t
  G(\oo_t,\vec\my,\mpp,\vec\mv,\mq)=\int_{\G_t}[\mathcal{W}_1(\vec\my,\vec\mv)
  +\mathcal{W}_2(\vec\my,\mpp,\vec\mv,\mq)]\langle\vec
  V,\vn_t\rangle\ds_t
\end{equation}
where
\begin{eqnarray*}
  \mathcal{W}_1(\vec\my,\vec\mv)&\defmath&\frac{1}{2}\seminorm{\vec\my-\vyd}^2+\alpha\md(\vec\my-\vg):\md\vec\mv;\\
 \mathcal{W}_2(\vec\my,\mpp,\vec\mv,\mq)&\defmath&\langle
 \alpha\Delta\vec\my-\n \mpp+\vf,\vec\mv\rangle+\langle\vec\my-\vg,\alpha\Delta\vec\mv-\n\mq\rangle+\mq\mdiv\vg.
\end{eqnarray*}
and $\vn_t$ denotes the outward unit normal to the boundary $\G_t$.

 By the previous choice of $\vf, \vg$ and $\vyd$, and $\vec V\in \mathcal{D}^1(\rn,\rn)$, $\p_t
G(\oo_t,\vec\my,\mpp,\vec\mv,\mq)$ exists everywhere in $[0,\tau]$
for all $(\vec\my,\mpp,\vec\mv,\mq)\in X\times Y.$ Hence (H2) is
satisfied.

For sufficiently smooth domains $\oo$ and vector fields $\vec
V\in\mathcal{D}^1(\rn,\rn)$, we have shown that $(\vy^t,p^t)$
(resp., $(\vvv^t,q^t)$) converge to $(\vy,p)$ (resp., $(\vvv,q)$) in
the $H^2\times H^1-$strong topology as $t$ goes to zero in the
previous section. Hence
\begin{equation*}
  \vec \my_t\rightarrow \vec\my=\Pi\, \vy,\;\mbox{ and }
   \; \vec \mv_t\rightarrow \vec\mv=\Pi\, \vvv\qquad\mbox{ strongly in }\;
   H^2(\rn)^N,
\end{equation*}
and
\begin{equation*}
  \mpp_t\rightarrow \mpp=\Pi \,p,\; \mbox{ and  }\mq_t\rightarrow\mq=\Pi q\;\qquad \mbox{strongly in}\;
  H^1(\rn).
\end{equation*}
by the following lemma.
\begin{lemma}[see\cite{delfour}]\label{fse:lemma}
 For any integer $m\geq 1$,  the velocity field $V\in \mathcal{D}^m(\rn,\rn)$ and a function $\Phi\in H^m(\rn),$
 if
 \begin{equation*}
   y^t\rightarrow y^0 \qquad \mbox{in}\;\; H^m(\oo)\mbox{-strong}
 \end{equation*}
 we have
 \begin{equation*}
   Y_t\rightarrow Y_0 \qquad \mbox{in}\;\; H^m(\rn)\mbox{-strong}
 \end{equation*}
 where $Y_t:=(\Pi\,y^t)\circ T_t^{-1}$. We also can show that the
 above result also holds for the weak topology of $H^m(\rn)$.
\end{lemma}
Furthermore, assumptions(H3)(i) and (H4)(i) are satisfied for the
$H^3\times H^2-$strong topology.

Now let's check (H3)(ii) and (H4)(ii).  Since
$(\vec\my,\mpp,\vec\mv,\mq)\in X\times Y$ and $\oo_t$ is
sufficiently smooth, we can use Stokes' formula to rewrite
(\ref{fse:aaa}) as:
\begin{equation*}
\p_t
G(\oo_t,\vec\my,\mpp,\vec\mv,\mq)=\int_{\oo_t}\mdiv\{[\mathcal{W}_1(\vec\my,\vec\mv)
  +\mathcal{W}_2(\vec\my,\mpp,\vec\mv,\mq)]\,\vec
V\}\dx,\;\forall\,(\vec\my,\mpp)\in X, (\vec\mv,\mq)\in Y.
\end{equation*}
Now introduce the mapping
\begin{equation*}
  (\vec\my,\mpp,\vec\mv,\mq)\mapsto [\mathcal{W}_1(\vec\my,\vec\mv)
  +\mathcal{W}_2(\vec\my,\mpp,\vec\mv,\mq)]\,\vec
V :\; X\times
  Y\rightarrow H^1(\rn)^N
\end{equation*}
which is linear and continuous.

Furthermore, by transformation $T_t$, the mapping
 \begin{equation*}
 (t;\vec\my,\mpp,\vec\mv,\mq)\mapsto\p_t
G(\oo_t,\vec\my,\mpp,\vec\mv,\mq)=\int_\oo
\mdiv\{[\mathcal{W}_1(\vec\my,\vec\mv)
  +\mathcal{W}_2(\vec\my,\mpp,\vec\mv,\mq)]\,\vec
V \}\circ T_t\,J_t\dx
 \end{equation*}
from $[0,\tau]\times X\times Y$ to $\mathbb{R}$ is continuous and
(H3)(ii) and (H4)(ii) are verified. This completes the verification
of the four assumptions of Theorem \ref{fsp:correa}.

Hence we obtain
\begin{equation}\label{dd}
  \diff J(\oo;\vec V)=\inf_{(\vec\my,\mpp)\in X(0)}\sup_{(\vec\mv,\mq)\in Y(0)}\p_t
  G(\oo_t,\vec\my,\mpp,\vec\mv,\mq)|_{t=0}.
\end{equation}
We also note that the expression (\ref{fse:aaa}) is a boundary
integral on $\G_t$ which will not depend on $(\vec\my,\mpp)$ and
$(\vec\mv,\mq)$ outside of $\overline{\oo}_t$, so the inf and the
sup in (\ref{dd}) can be dropped, we then get
\begin{equation*}
  \diff J(\oo;\vec V)=\int_\G [\mathcal{W}_1(\vy,\vvv)+\mathcal{W}_2(\vy,p,\vvv,q)]\langle\vec V,\vn\rangle\ds
\end{equation*}
However, $\vy=\vg$, $\vp=0$ and (\ref{comp}) imply
$\mathcal{W}_2(\vy,p,\vvv,q)=0$ on the boundary $\G$. Finally we
have
\begin{equation*}
  \diff J(\oo;\vec V)=\int_\G\left\{\frac{1}{2}\seminorm{\vy-\vyd}^2
  +\alpha\md(\vy-\vg):\md\vvv\right\}\vv\cdot\vn\ds
\end{equation*}
We also find that the expression of Eulerian derivative obtained by
function space embedding  was the same as (\ref{end}) which was
obtained by the function space parametrization technique, but the
second method is obviously quick.
\section{Gradient algorithm and numerical implementation}
In this section, we will give a gradient type algorithm and some
numerical examples in two dimensions to prove that our previous
methods (i.e. Function Space Parametrization \& Function Space
Embedding) could be very useful and efficient for the numerical
implementation of shape problems.

We describe a gradient type algorithm for the minimization of a cost
function $J(\oo)$. As we have just seen, the general form of its
Eulerian derivative is
\begin{equation*}
  \diff J(\oo;\vec V)=\int_\G v \vec V\cdot\vn\ds
\end{equation*}
where $v$ is given by a result like (\ref{end}). Ignoring
regularization, a descent direction is found by defining
\begin{equation}
  \vec V=-v\;\vn
\end{equation}
and then we can update the shape $\oo$ as
\begin{equation}
  \oo_k=(\mathrm{Id}+h_k\vec V)\oo
\end{equation}
where $h_k$ is a small descent step at $k$-th iteration.

There are also other choices for the definition of the descent
direction. The method used in this paper is to change the scalar
product with respect to which we compute a descent direction, for
instance, $H^1(\oo)^2$. In this case, the descent direction is the
unique element $\vec d\in H^1(\oo)^2$ such that for every $\vec V\in
H^1(\oo)^2,$
\begin{equation}\label{reg}
  \int_\oo\md\vec d :\md \vec V\dx=\diff J(\oo;\vec V).
\end{equation}
  The computation of $\vec d$ can also be interpreted as a regularization
  of the shape gradient, and the choice of $H^1(\oo)^2$ as space of
  variations is more dictated by technical considerations rather
  than theoretical ones.

The resulting algorithm can be summarized as follows:
\begin{itemize}
    \item [(1)] Choose an initial shape $\oo_0$;
    \item [(2)] Compute the state system and adjoint state system, then
    we can evaluate the descent direction $\vec d_k$ by using (\ref{reg})
    with $\oo=\oo_k;$
    \item[(3)] Set $\oo_{k+1}=(\mathrm{Id}-h_k\vec d_k) \,\oo_k,$ where $h_k$
    is a small positive real number and can be chosen by some rules, such as Armijo rule.
\end{itemize}

Our numerical solutions are obtained under FreeFem++ \cite{hecht}.
To illustrate the theory, we have solved the following minimization
problem
\begin{equation}\label{exam:fun}
\min_{\oo}\frac{1}{2}\int_{\oo}(\vy-\vyd)^2\dx
\end{equation}
subject to
\begin{equation}
  \left\{%
  \begin{array}{ll}
-\alpha\Delta\vy+\n p=\vf&\quad\mbox{in}\;\oo\\
\mdiv \vy=0&\quad\mbox{in}\;\oo\\
\vy=0&\quad\mbox{on}\;\G;\\
  \end{array}%
  \right.
\end{equation}
The domain $\oo$ is an annuli, and its boundary has two part: the
outer boundary $\G_{out}$ is a unit circle which is fixed; the inner
boundary $\G_{in}$ which is to be optimized. We choose the target
velocity $\vyd=(\vy_{1\diff},\vy_{2\diff})$ as follows:
$$\vy_{1\diff}=-\frac{y(\sqrt{x^2+y^2}-0.2)(\sqrt{x^2+y^2}-1)}{\sqrt{x^2+y^2}},
\;\vy_{2\diff}=\frac{x(\sqrt{x^2+y^2}-0.2)(\sqrt{x^2+y^2}-1)}{\sqrt{x^2+y^2}},$$
and the target inner boundary $\G_{in}$ is a concentric circle with
radius $0.2$.
We will solve this model problem with two different initial shapes:\\
\noindent{\bf Case 1:} A circle whose center is at origin with
radius 0.4, i.e., $x^2+y^2=0.4^2$;\\
\noindent{\bf Case 2:} A parabolic: ${x^2}/{9}+{y^2}/{4}={1}/{25}$.

Now the initial mesh of the two cases are shown in \autoref{fig1}
and \autoref{fig2}.

\begin{figure}[!htbp]
\renewcommand{\captionlabelfont}{\small}
\setcaptionwidth{2in}
\begin{minipage}[b]{0.50\textwidth}
  \centering
  \includegraphics[width=2.5in]{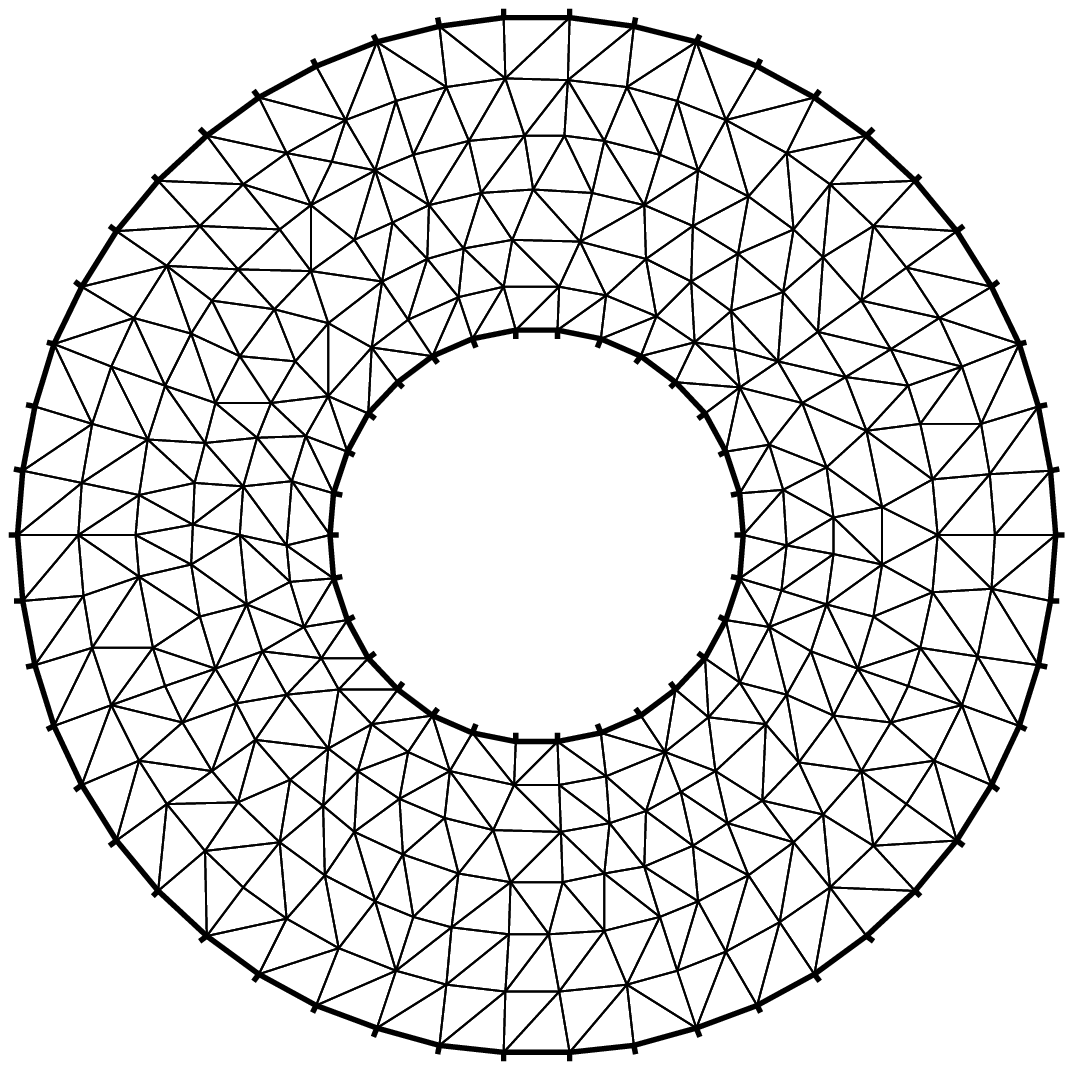}
  \caption{Initial mesh in Case 1 with 292 nodes.\label{fig1}}
  \end{minipage}
\begin{minipage}[b]{0.50\textwidth}
  \centering
  \includegraphics[width=2.5in]{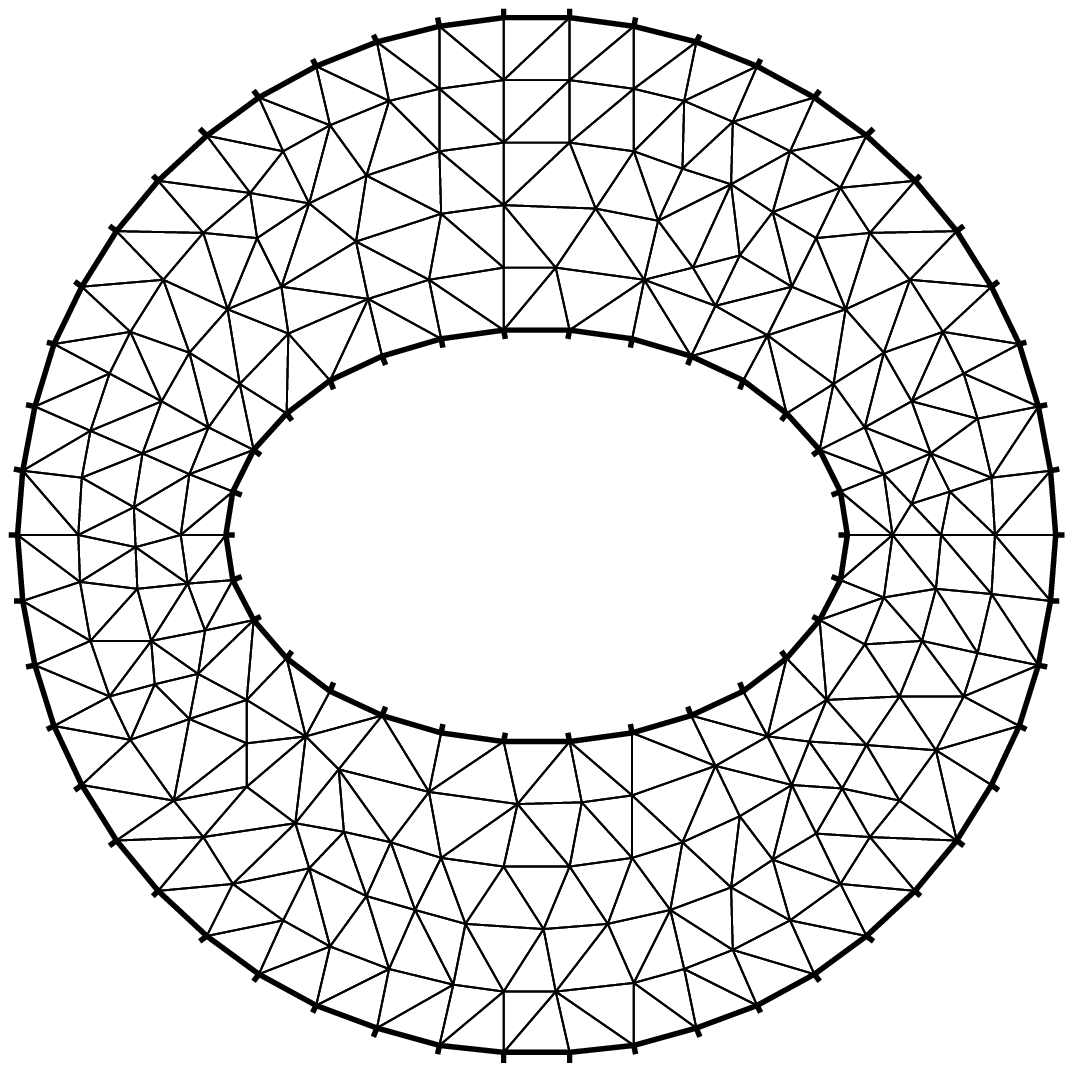}
  \caption{Initial mesh in Case 2 with 226 nodes.\label{fig2}}
  \end{minipage}
\end{figure}

 We will use the mixed finite element method to solve
the state system (\ref{p}) and adjoint state system (\ref{u}) on a
triangular mesh, and the popular P1-bubble/P1 finite element couple
(see \cite{girault86}) is chosen for the velocity-pressure couple.
We run the program on a home PC.

\begin{figure}[!htbp]
\centering
  \includegraphics[width=0.86\textwidth]{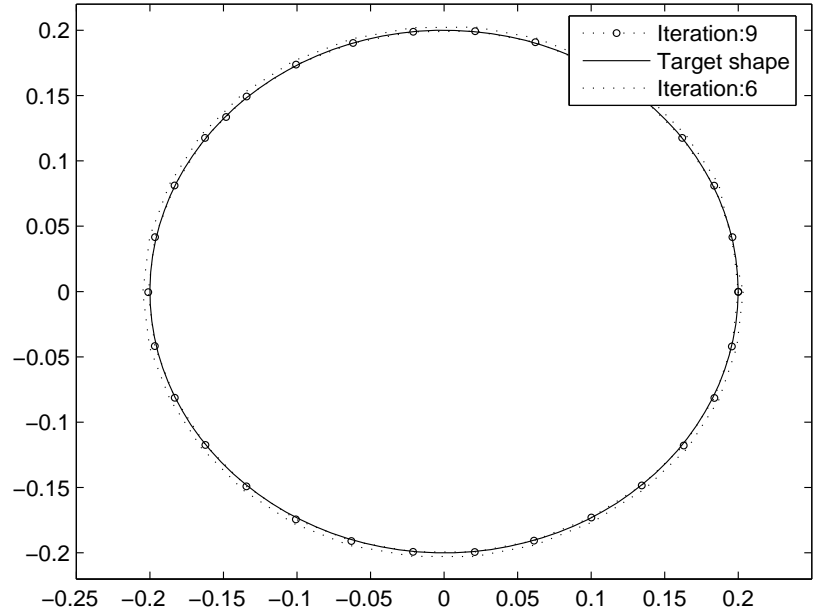}
  \caption{Case 1, $\alpha=1$, CPU time after 9 iterations: 37.235 second\label{yuan:fig1}}
\end{figure}
\begin{figure}[!htbp]
\centering
  \includegraphics[width=0.86\textwidth]{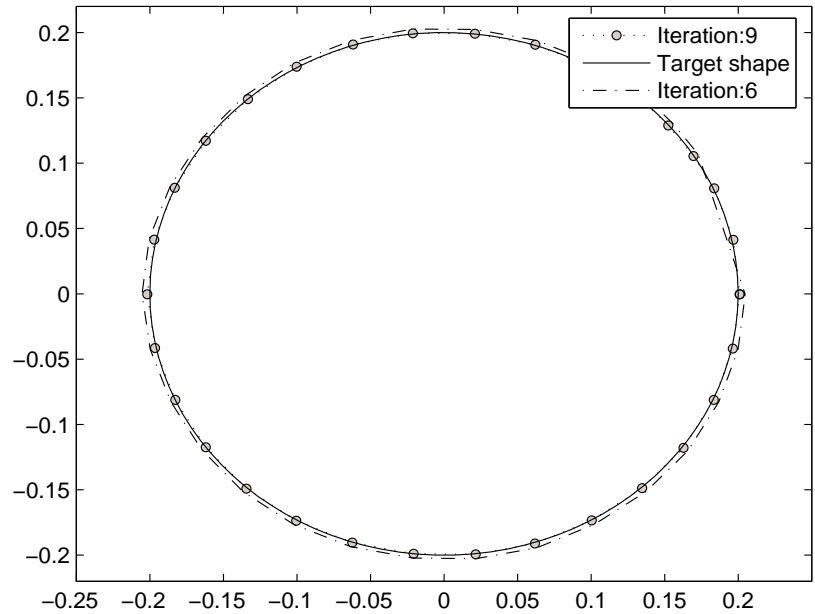}
  \caption{Case 1, $\alpha=0.1$, CPU time after 9 iterations: 36.125 second\label{yuan:fig2}}
\end{figure}

\begin{figure}[!htbp]
\centering
  \includegraphics[width=0.8\textwidth]{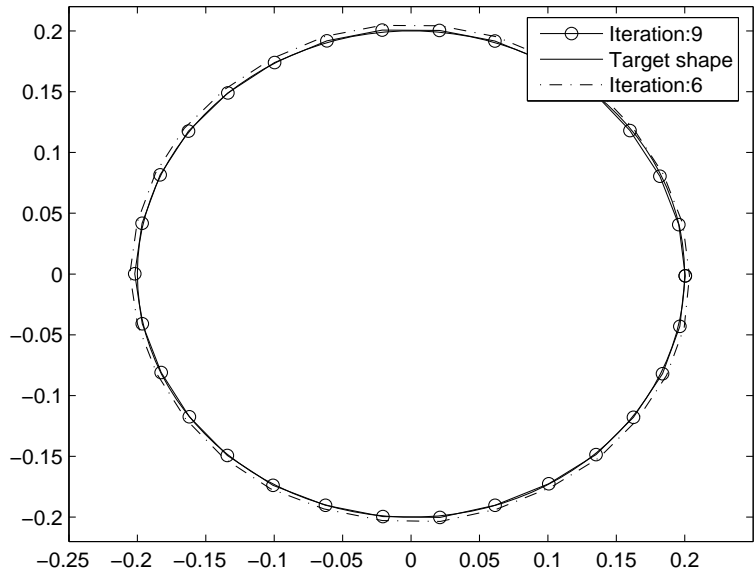}
  \caption{Case 1, $\alpha=0.01$, CPU time after 9 iterations: 37.14 second\label{yuan:fig3}}
\end{figure}

\begin{figure}[!htbp]
\centering
  \includegraphics[width=0.8\textwidth]{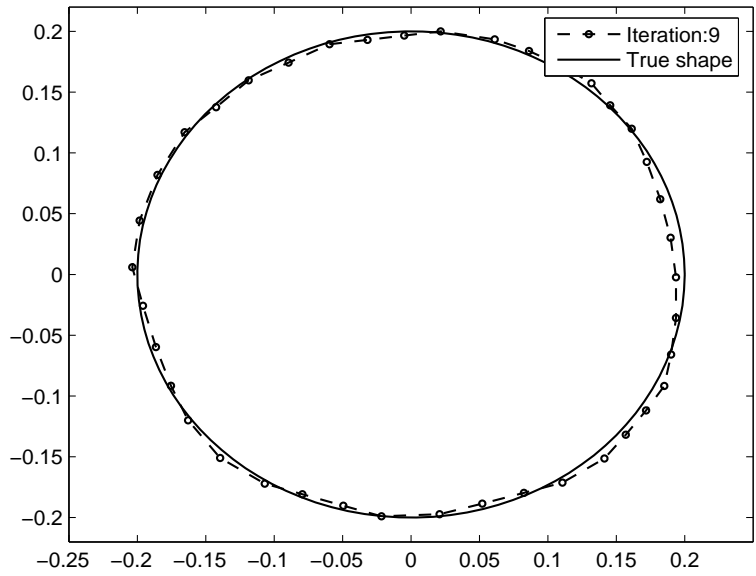}
  \caption{Case 1, $\alpha=0.001$, CPU time after 9 iterations: 47.375 second\label{yuan:fig4}}
\end{figure}
In Case 1, \autoref{yuan:fig1}---\autoref{yuan:fig4} give the
comparison between the target shape with iterated shape for the
viscosity coefficient $\alpha=1, 0.1, 0.01, 0.001$, respectively. We
can find that for $\alpha=1, 0.1, 0.01,$ we have nice
reconstruction, but for $\alpha=0.001,$ the result is not so
satisfied in \autoref{yuan:fig4}.

\begin{figure}[!htbp]
\centering
  \includegraphics[width=0.82\textwidth]{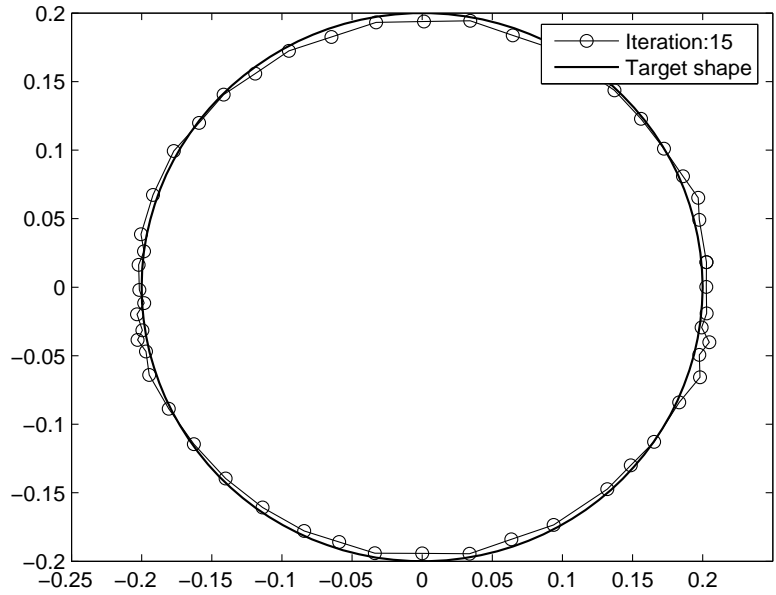}
  \caption{Case 2, $\alpha=1$, CPU time after 15 iterations: 64.11 second\label{tuoyuan:fig1}}
\end{figure}
\begin{figure}[!htbp]
\centering
  \includegraphics[width=0.82\textwidth]{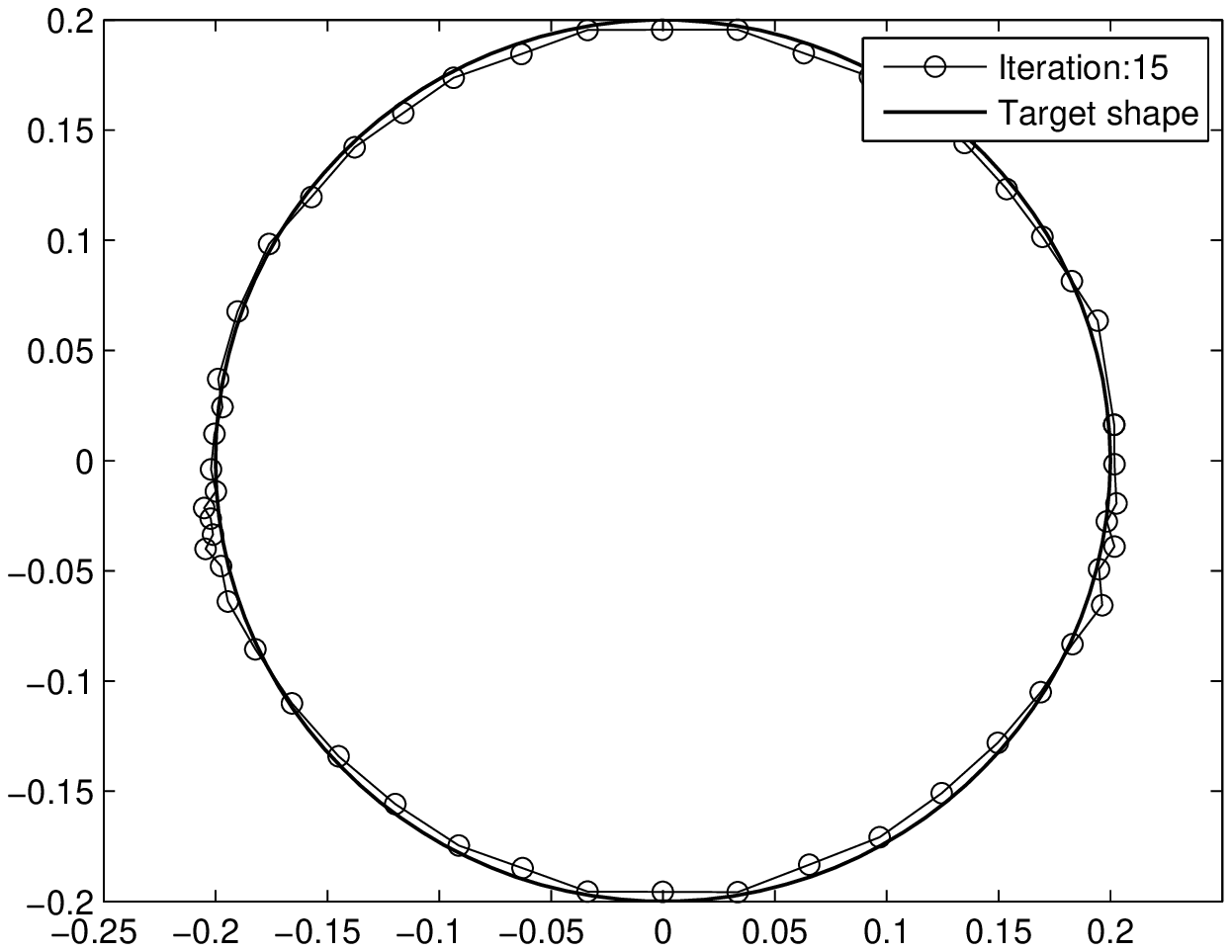}
  \caption{Case 2, $\alpha=0.1$, CPU time after 15 iterations: 64.172 second\label{tuoyuan:fig2}}
\end{figure}
\begin{figure}[!htbp]
\centering
  \includegraphics[width=0.8\textwidth]{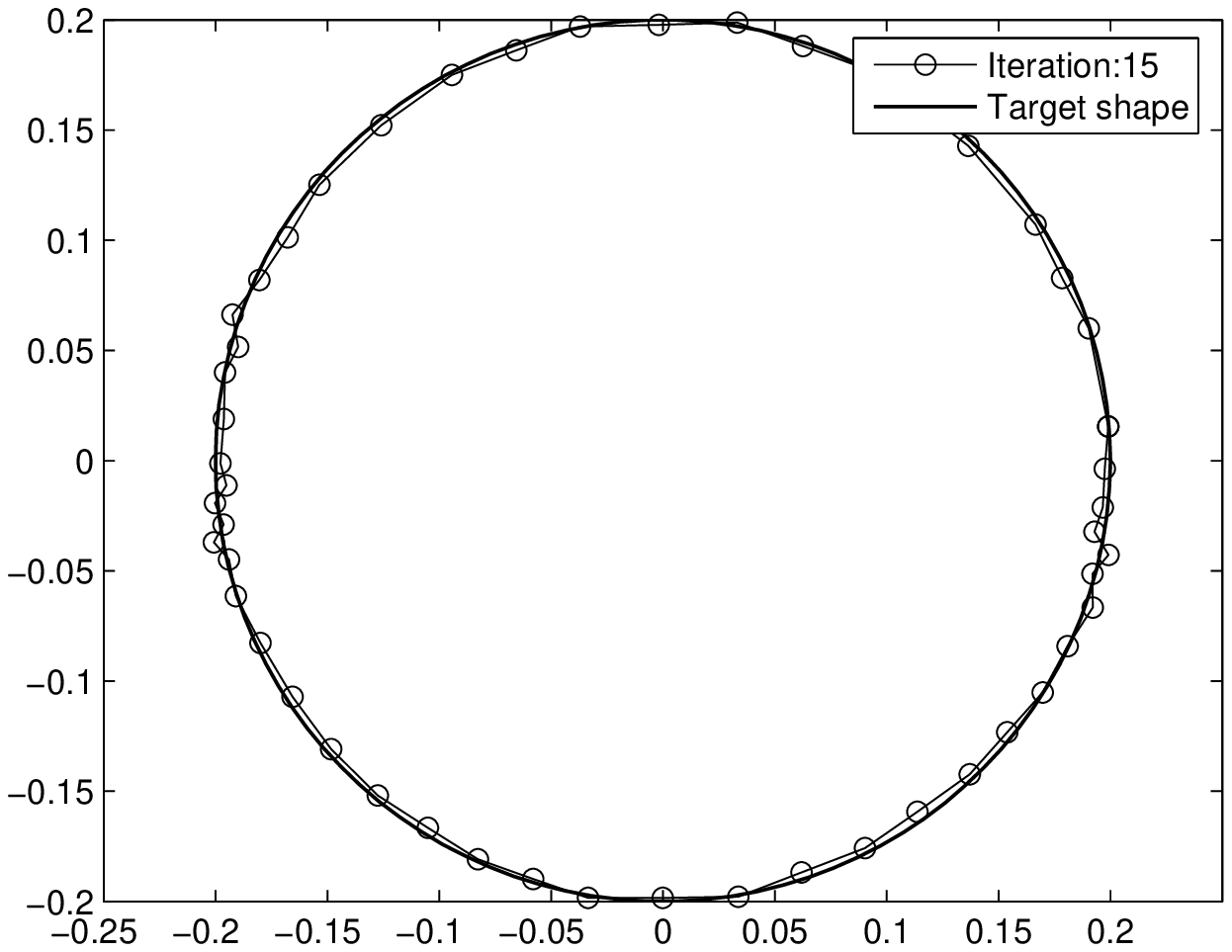}
  \caption{Case 2, $\alpha=0.01$, CPU time after 15 iterations: 66.141 second\label{tuoyuan:fig3}}
\end{figure}
\begin{figure}[!htbp]
\centering
  \includegraphics[width=0.8\textwidth]{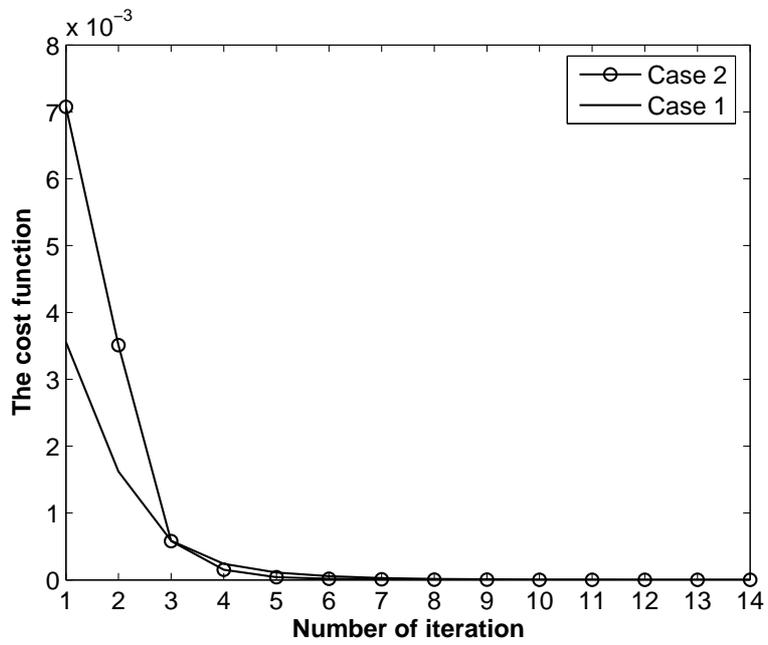}
  \caption{Convergence history of the cost function in two cases with $\alpha=0.01$.\label{fig5}}
\end{figure}
 In Case 2,
\autoref{tuoyuan:fig1}---\autoref{tuoyuan:fig3} represent the
comparison between the target shape with iterated shape for the
viscosity coefficient $\alpha=1, 0.1, 0.01$, respectively. It can be
shown that for fixed viscosity, Case 1 has better reconstruction
than Case 2, that's to say, the iteration process depends on the
choice of the initial shape.

 \autoref{fig5} shows the fast
convergence of our cost function (\ref{exam:fun}) in Case 1 and Case
2 for the viscosity $\alpha=0.01$.

Finally, the numerical examples show the feasibility of the proposed
iteration algorithm and further research is necessary on efficient
implementations.

\end{document}